\documentclass[a4paper, 12 pt]{article}

\usepackage[T2A]{fontenc}
\usepackage[cp1251]{inputenc}
\usepackage[english]{babel}
\usepackage[tbtags]{amsmath}
\usepackage{amsfonts,amssymb}

\usepackage{mathrsfs}

\usepackage{geometry} % Меняем поля страницы
\begin{document}
\title{Generalized Riordan arrays and zero generalized Pascal matrices}
 
 \author{E. Burlachenko}
 \date{}
 
 \maketitle

\begin{abstract}

Generalized Pascal matrix whose elements are generalized binomial coefficients is included in the group of  generalized Riordan arrays. There is a special set of generalized Riordan arrays defined by parameter $q$. If $q=0$, they are ordinary Riordan arrays, if $q=1$, they are exponential Riordan arrays. In other cases, except $q=-1$, they are arrays associated with the $q$-binomial coefficients as well as the exponential Riordan arrays are associated with the ordinary binomial coefficients. Case $q=-1$ does not fit into the concept of generalized Riordan arrays, but it is necessary to expand for it. Introduced a special class of matrices, each of which is a limiting case of a certain set of generalized Pascal matrices. It is shown that every such matrix included in the matrix group similar to the generalized Riordan group.

\end{abstract}
\section{Introduction}

Transformations, corresponding to multiplication and composition of  series,
play the main role in the space of formal power series over the field of real or complex  numbers. Multiplication is geven by the matrix $\left( a\left( x \right),x \right)$ $n$th column of which, $n=0,\text{ }1,\text{ }2,\text{ }...$ ,  has the generating function $b\left( x \right){{x}^{n}}$; composition is given by the matrix $\left( 1,a\left( x \right) \right)$ $n$th column of which  has the generating function ${{a}^{n}}\left( x \right)$, ${{a}_{0}}=0$:
$$\left( b\left( x \right),x \right)g\left( x \right)=b\left( x \right)g\left( x \right), \qquad\left( 1,a\left( x \right) \right)g\left( x \right)=g\left( a\left( x \right) \right).$$
Matrix
$$\left( b\left( x \right),x \right)\left( 1,a\left( x \right) \right)=\left( b\left( x \right),a\left( x \right) \right)$$
is called Riordan array [1] – [4]; $n$th column of Riordan array has the generating function $b\left( x \right){{a}^{n}}\left( x \right)$. Thus,
$$\left( b\left( x \right),a\left( x \right) \right)f\left( x \right){{g}^{n}}\left( x \right)=b\left( x \right)f\left( a\left( x \right) \right){{g}^{n}}\left( a\left( x \right) \right),$$
$$\left( b\left( x \right),a\left( x \right) \right)\left( f\left( x \right),g\left( x \right) \right)=\left( b\left( x \right)f\left( a\left( x \right) \right),g\left( a\left( x \right) \right) \right).$$
Matrices $\left( b\left( x \right),a\left( x \right) \right)$, ${{b}_{0}}\ne 0$, ${{a}_{1}}\ne 0$, form a group called the Riordan group.

$n$th coefficient of the series $a\left( x \right)$, $\left( n,m \right)$th element of the matrix $A$, $n$th row and $n$th column of the matrix $A$ will be denoted  respectively by 
$$\left[ {{x}^{n}} \right]a\left( x \right), \qquad {{\left( A \right)}_{n,m}}, \qquad  \left[ n,\to  \right]A,  \qquad \left[\uparrow ,n \right]A.$$
We associate rows and columns of matrices with the generating functions of their elements. For the elements of the lower triangular matrices will be appreciated that ${{\left( A \right)}_{n,m}}=0$, if $n<m$. 

Matrices 
$${{\left| {{e}^{x}} \right|}^{-1}}\left( b\left( x \right),a\left( x \right) \right)\left| {{e}^{x}} \right|={{\left( b\left( x \right),a\left( x \right) \right)}_{{{e}^{x}}}},$$
where $\left| {{e}^{x}} \right|$ is the diagonal matrix whose diagonal elements are equal to the coefficients of  the series ${{e}^{x}}$: $\left| {{e}^{x}} \right|a\left( x \right)=\sum\nolimits_{n=0}^{\infty }{{{{a}_{n}}{{x}^{n}}}/{n!}\;}$, are called exponential Riordan arrays. Denote
$$\left[ n,\to  \right]{{\left( b\left( x \right),a\left( x \right) \right)}_{{{e}^{x}}}}={{s}_{n}}\left( x \right), \qquad{{b}_{0}}\ne 0, \qquad{{a}_{1}}\ne 0.$$
Then
$${{\left( b\left( x \right),a\left( x \right) \right)}_{{{e}^{x}}}}{{\left( 1-\varphi x \right)}^{-1}}={{\left| {{e}^{x}} \right|}^{-1}}\left( b\left( x \right),a\left( x \right) \right){{e}^{\varphi x}}={{\left| {{e}^{x}} \right|}^{-1}}b\left( x \right)\exp \left( \varphi a\left( x \right) \right),$$
or
$$\sum\limits_{n=0}^{\infty }{\frac{{{s}_{n}}\left( \varphi  \right)}{n!}{{x}^{n}}}=b\left( x \right)\exp \left( \varphi a\left( x \right) \right).$$
Sequence of polynomials ${{s}_{n}}\left( x \right)$ is called Sheffer sequence and is a subject of study of the classical umbral calculus [5]. Examples of Sheffer polynomials are  the Bernoulli, Euler, Hermite, Laguerre polynomials. They correspond to the matrices
$${{\left( \frac{x}{{{e}^{x}}-1},x \right)}_{{{e}^{x}}}}, \quad{{\left( \frac{2}{{{e}^{x}}+1},x \right)}_{{{e}^{x}}}}, \quad{{\left( {{e}^{-{{x}^{2}}}},2x \right)}_{{{e}^{x}}}}, \quad{{\left( \frac{1}{1-x},\frac{-x}{1-x} \right)}_{{{e}^{x}}}}.$$

Matrix
$$P=\left( \frac{1}{1-x},\frac{x}{1-x} \right)={{\left( {{e}^{x}},x \right)}_{{{e}^{x}}}}=\left( \begin{matrix}
   1 & 0 & 0 & 0 & 0 & \ldots   \\
   1 & 1 & 0 & 0 & 0 & \ldots   \\
   1 & 2 & 1 & 0 & 0 & \ldots   \\
   1 & 3 & 3 & 1 & 0 & \ldots   \\
   1 & 4 & 6 & 4 & 1 & \ldots   \\
   \vdots  & \vdots  & \vdots  & \vdots  & \vdots  & \ddots   \\
\end{matrix} \right).$$
has a special status and is called Pascal matrix. Power of the Pascal matrix is defined by the identity

$${{P}^{\varphi }}=\left( \frac{1}{1-\varphi x},\frac{x}{1-\varphi x} \right)={{\left( {{e}^{\varphi x}},x \right)}_{{{e}^{x}}}}.$$

Matrices 
$${{\left| c\left( x \right) \right|}^{-1}}\left( b\left( x \right),a\left( x \right) \right)\left| c\left( x \right) \right|={{\left( b\left( x \right),a\left( x \right) \right)}_{c\left( x \right)}},$$
where $\left| c\left( x \right) \right|$ is the diagonal matrix: $\left| c\left( x \right) \right|a\left( x \right)=\sum\nolimits_{n=0}^{\infty }{{{c}_{n}}{{a}_{n}}}{{x}^{n}}$, ${{c}_{n}}\ne 0$, are called generalized Riordan arrays [3]  ($\left( c \right)$-Riordan arrays[4]). Denote
$$\left[ n,\to  \right]{{\left( b\left( x \right),a\left( x \right) \right)}_{c\left( x \right)}}={{u}_{n}}\left( x \right), \qquad{{b}_{0}}\ne 0, \qquad{{a}_{1}}\ne 0.$$
Then
$$\sum\limits_{n=0}^{\infty }{{{c}_{n}}{{u}_{n}}}\left( \varphi  \right){{x}^{n}}=b\left( x \right)c\left( \varphi a\left( x \right) \right).$$
Sequence of polynomials ${{u}_{n}}\left( x \right)$ is called Boas-Buck sequence and is a subject of study of the non-classical umbral calculus [5], [6]. Examples of Boas-Buck polynomials are  the Chebyshev polynomials of the first and second kind. They correspond to the ordinary Riordan arrays ($c\left( x \right)={{\left( 1-x \right)}^{-1}}$)
$$\frac{1}{2}\left( \frac{1-{{x}^{2}}}{1+{{x}^{2}}},\frac{2x}{1+{{x}^{2}}} \right), \qquad\left( \frac{1}{1+{{x}^{2}}},\frac{2x}{1+{{x}^{2}}} \right),$$
where it is believed
$${{T}_{0}}\left( x \right)=\frac{1}{2},  \qquad 2\sum\limits_{n=0}^{\infty }{{{T}_{n}}\left( \varphi  \right){{x}^{n}}}=\frac{1-{{x}^{2}}}{1-2\varphi x+{{x}^{2}}}.$$
Laguerre polynomials give an example of the ambiguous Boas-Buck structure [6].  Corresponding matrix can be represented as anexponential Riordan array, and as
$${{\left( {{e}^{x}},-x \right)}_{c\left( x \right)}}, \qquad c\left( x \right)=\sum\limits_{n=0}^{\infty }{\frac{{{x}^{n}}}{{{\left( n! \right)}^{2}}}}.$$

Generalized Riordan arrays are associated with the following generalization of the binomial coefficients [7]. For the coefficients of the formal power series  $b\left( x \right)$, ${{b}_{0}}=0$; ${{b}_{n}}\ne 0$,  $n>0$, denote
$${{b}_{0}}!=1,\qquad {{b}_{n}}!=\prod\limits_{m=1}^{n}{{{b}_{m}}},\qquad {{\left( \begin{matrix}
   n  \\
   m  \\
\end{matrix} \right)}_{b}}=\frac{{{b}_{n}}!}{{{b}_{m}}!{{b}_{n-m}}!};\qquad   {{\left( \begin{matrix}
   n  \\
   m  \\
\end{matrix} \right)}_{b}}=0,  \qquad m>n.$$
Then
\[{{\left( \begin{matrix}
   n  \\
   m  \\
\end{matrix} \right)}_{b}}={{\left( \begin{matrix}
   n-1  \\
   m-1  \\
\end{matrix} \right)}_{b}}+\frac{{{b}_{n}}-{{b}_{m}}}{{{b}_{n-m}}}{{\left( \begin{matrix}
   n-1  \\
   m  \\
\end{matrix} \right)}_{b}}.\]
Consider matrix
$${{P}_{c\left( x \right)}}={{\left( c\left( x \right),x \right)}_{c\left( x \right)}}= \left( \begin{matrix}
   \frac{{{c}_{0}}{{c}_{0}}}{{{c}_{0}}} & 0 & 0 & 0 & 0 & \ldots   \\
   \frac{{{c}_{0}}{{c}_{1}}}{{{c}_{1}}} & \frac{{{c}_{1}}{{c}_{0}}}{{{c}_{1}}} & 0 & 0 & 0 & \ldots   \\
   \frac{{{c}_{0}}{{c}_{2}}}{{{c}_{2}}} & \frac{{{c}_{1}}{{c}_{1}}}{{{c}_{2}}} & \frac{{{c}_{2}}{{c}_{0}}}{{{c}_{2}}} & 0 & 0 & \ldots   \\
   \frac{{{c}_{0}}{{c}_{3}}}{{{c}_{3}}} & \frac{{{c}_{1}}{{c}_{2}}}{{{c}_{3}}} & \frac{{{c}_{2}}{{c}_{1}}}{{{c}_{3}}} & \frac{{{c}_{3}}{{c}_{0}}}{{{c}_{3}}} & 0 & \ldots   \\
   \frac{{{c}_{0}}{{c}_{4}}}{{{c}_{4}}} & \frac{{{c}_{1}}{{c}_{3}}}{{{c}_{4}}} & \frac{{{c}_{2}}{{c}_{2}}}{{{c}_{4}}} & \frac{{{c}_{3}}{{c}_{1}}}{{{c}_{4}}} & \frac{{{c}_{4}}{{c}_{0}}}{{{c}_{4}}} & \ldots   \\
   \vdots  & \vdots  & \vdots  & \vdots  & \vdots  & \ddots   \\
\end{matrix} \right),\quad
{{\left( {{P}_{c\left( x \right)}} \right)}_{n,m}}=\frac{{{c}_{m}}{{c}_{n-m}}}{{{c}_{n}}}.$$  
 Denote $[\uparrow ,1]{{P}_{c\left( x \right)}}=b\left( x \right)$. If ${{c}_{0}}=1$,  then
$${{c}_{n}}=\frac{c_{1}^{n}}{{{b}_{n}}!} ,\qquad {{\left( {{P}_{c\left( x \right)}} \right)}_{n,m}}={{\left( \begin{matrix}
   n  \\
   m  \\
\end{matrix} \right)}_{b}}.$$
Let ${{c}_{n}}\in \mathbb{R}$, ${{c}_{0}}=1$. Since ${{P}_{c\left( x \right)}}={{P}_{c\left( \varphi x \right)}}$, we take for uniqueness that ${{c}_{1}}=1$. Matrix ${{P}_{c\left( x \right)}}$ will be called generalized Pascal matrix.

In Section 2 we consider the set of generalized Pascal matrices as a group under  Hadamard multiplication and introduce a special system of matrices, which implies the concept of zero generalized Pascal matrices. In Section 3 we will give an idea of the algebra associated with such matrix. In Sections 4-6 we consider the main varieties of these matrices, in particular, fractal zero generalized Pascal matrices, an example of which is the Pascal triangle modulo 2. Emphasis is on the algebras associated with these matrices. In Section 7 we will give an idea of the matrix group similar generalized  Riordan group which includes zero generalized Pascal matrix.

\section{Special system of generalized Pascal matrices}

Elements of the matrix ${{P}_{c\left( x \right)}}$, – denote them ${{\left( {{P}_{c\left( x \right)}} \right)}_{n,m}}={{{c}_{m}}{{c}_{n-m}}}/{{{c}_{n}}}\;=\left( n,m \right)$ for generality which will be discussed later, –  satisfy the identities 
	$$\left( n,0 \right)=1, \qquad\left( n,m \right)=\left( n,n-m \right), \eqno  	(1)$$
$$\text{  }\left( n+q,q \right)\left( n+p,m+p \right)\left( m+p,p \right)=\left( n+p,p \right)\left( n+q,m+q \right)\left( m+q,q \right),\eqno(2)$$
$q$, $p=0$, $1$, $2$, … It means that each matrix ${{P}_{c\left( x \right)}}$ can be associated with the algebra of formal power series whose elements are multiplied by the rule
$$a\left( x \right)\circ b\left( x \right)=g\left( x \right),  \qquad{{g}_{n}}=\sum\limits_{m=0}^{n}{\left( n,m \right){{a}_{m}}}{{b}_{m-n}},$$
that is, if ${{\left( A \right)}_{n,m}}={{a}_{n-m}}\left( n,m \right)$, ${{\left( B \right)}_{n,m}}={{b}_{n-m}}\left( n,m \right)$, ${{\left( G \right)}_{n,m}}={{g}_{n-m}}\left( n,m \right)$,   then $AB=BA=G$:
$${{g}_{n}}={{\left( n+p,p \right)}^{-1}}\sum\limits_{m=0}^{n}{\left( n+p,m+p \right)\left( m+p,p \right){{a}_{n-m}}{{b}_{m}}}=$$
$$={{\left( n+q,q \right)}^{-1}}\sum\limits_{m=0}^{n}{\left( n+q,m+q \right)\left( m+q,q \right){{a}_{n-m}}{{b}_{m}}}.$$

The set of generalized Pascal matrix is a group under Hadamard multiplication (we denote this operation $\times $):
$${{P}_{c\left( x \right)}}\times {{P}_{g\left( x \right)}}={{P}_{c\left( x \right)\times g\left( x \right)}},  \qquad c\left( x \right)\times g\left( x \right)=\sum\limits_{n=0}^{\infty }{{{c}_{n}}{{g}_{n}}}{{x}^{n}}.$$

Introduce the special system of matrices
$${}_{\varphi ,q}P={}_{q}P\left( \varphi  \right)={{P}_{c\left( \varphi ,q,x \right)}},  \qquad c\left( \varphi ,q,x \right)=\left( \sum\limits_{n=0}^{q-1}{{{x}^{n}}} \right){{\left( 1-\frac{{{x}^{q}}}{\varphi } \right)}^{-1}}, \qquad q>1.$$
Then 
$${{c}_{qn+i}}=\frac{1}{{{\varphi }^{n}}},\qquad  0\le i<q; \qquad {{c}_{qn-i}}=\frac{1}{{{\varphi }^{n-1}}},\qquad  0<i\le q,$$
$$\frac{{{c}_{qm+j}}{{c}_{q\left( n-m \right)+i-j}}}{{{c}_{qn+i}}}=\frac{{{\varphi }^{n}}}{{{\varphi }^{m}}{{\varphi }^{n-m}}}=1,\qquad i\ge j;  \qquad=\frac{{{\varphi }^{n}}}{{{\varphi }^{m}}{{\varphi }^{n-m-1}}}=\varphi , \qquad i<j,$$
or
$${{\left( _{\varphi ,q}P \right)}_{n,m}}=1, \qquad n\left( \bmod q \right)\ge m\left( \bmod q \right); \qquad=\varphi , \qquad n\left( \bmod q \right)<m\left( \bmod q \right).$$
For example, $_{\varphi ,2}P$, $_{\varphi ,3}P$:
$$\left( \begin{matrix}
   1 & 0 & 0 & 0 & 0 & 0 & 0 & 0 & 0 & \ldots   \\
   1 & 1 & 0 & 0 & 0 & 0 & 0 & 0 & 0 & \ldots   \\
   1 & \varphi  & 1 & 0 & 0 & 0 & 0 & 0 & 0 & \ldots   \\
   1 & 1 & 1 & 1 & 0 & 0 & 0 & 0 & 0 & \ldots   \\
   1 & \varphi  & 1 & \varphi  & 1 & 0 & 0 & 0 & 0 & \ldots   \\
   1 & 1 & 1 & 1 & 1 & 1 & 0 & 0 & 0 & \dots   \\
   1 & \varphi  & 1 & \varphi  & 1 & \varphi  & 1 & 0 & 0 & \ldots   \\
   1 & 1 & 1 & 1 & 1 & 1 & 1 & 1 & 0 & \ldots   \\
   1 & \varphi  & 1 & \varphi  & 1 & \varphi  & 1 & \varphi  & 1 & \ldots   \\
   \vdots  & \vdots  & \vdots  & \vdots  & \vdots  & \vdots  & \vdots  & \vdots  & \vdots  & \ddots   \\
\end{matrix} \right), \qquad \left( \begin{matrix}
   1 & 0 & 0 & 0 & 0 & 0 & 0 & 0 & 0 & \ldots   \\
   1 & 1 & 0 & 0 & 0 & 0 & 0 & 0 & 0 & \ldots   \\
   1 & 1 & 1 & 0 & 0 & 0 & 0 & 0 & 0 & \ldots   \\
   1 & \varphi  & \varphi  & 1 & 0 & 0 & 0 & 0 & 0 & \ldots   \\
   1 & 1 & \varphi  & 1 & 1 & 0 & 0 & 0 & 0 & \ldots   \\
   1 & 1 & 1 & 1 & 1 & 1 & 0 & 0 & 0 & \ldots   \\
   1 & \varphi  & \varphi  & 1 & \varphi  & \varphi  & 1 & 0 & 0 & \ldots   \\
   1 & 1 & \varphi  & 1 & 1 & \varphi  & 1 & 1 & 0 & \ldots   \\
   1 & 1 & 1 & 1 & 1 & 1 & 1 & 1 & 1 & \ldots   \\
   \vdots  & \vdots  & \vdots  & \vdots  & \vdots  & \vdots  & \vdots  & \vdots  & \vdots  & \ddots   \\
\end{matrix} \right).$$

Elements of the matrix $_{\varphi ,q}P\times {{P}_{c\left( x \right)}}$ satisfy the identities (1), (2) for any values $\varphi $, so it makes sense to consider also the case $\varphi =0$ since it corresponds to a certain algebra of formal power series. 
It is clear that in this case the series $c\left( \varphi ,q,x \right)$ is not defined. Matrix $_{0,q}P\times {{P}_{c\left( x \right)}}$ and Hadamard product of such matrices will be called zero generalized Pascal matrix.\\
{\bfseries Remark.} Zero generalized Pascal matrix appears when considering the set of generalized Pascal matrices ${{P}_{g\left( q,x \right)}}$:
$${{\left( {{P}_{g\left( q,x \right)}} \right)}_{n,1}}=\left[ {{x}^{n}} \right]\frac{x}{\left( 1-x \right)\left( 1-qx \right)}=\sum\limits_{m=0}^{n-1}{{{q}^{m}}}, \qquad q\in \mathbb{R}.$$
Here $g\left( 0,x \right)={{\left( 1-x \right)}^{-1}}$,  $g\left( 1,x \right)={{e}^{x}}$. In other cases  (the $q$-umbral calculus [5]), except $q=-1$,
$$g\left( q,x \right)=\sum\limits_{n=0}^{\infty }{\frac{{{\left( q-1 \right)}^{n}}}{\left( {{q}^{n}}-1 \right)!}}{{x}^{n}},  \qquad\left( {{q}^{n}}-1 \right)!=\prod\limits_{m=1}^{n}{\left( {{q}^{m}}-1 \right)},  \qquad\left( {{q}^{0}}-1 \right)!=1.$$
Matrices  ${{P}_{g\left( q,x \right)}}$, $P_{g\left( q,x \right)}^{-1}$ also can  be defined as follows:
$$[\uparrow ,n]{{P}_{g\left( q,x \right)}}={{x}^{n}}\prod\limits_{m=0}^{n}{{{\left( 1-{{q}^{m}}x \right)}^{-1}}},    \qquad[n,\to ]P_{g\left( q,x \right)}^{-1}=\prod\limits_{m=0}^{n-1}{\left( x-{{q}^{m}} \right)}.$$
When  $q=-1$ we get the matrices ${{P}_{g\left( -1,x \right)}}$, $P_{g\left( -1,x \right)}^{-1}$:
$$\left( \begin{matrix}
   1 & 0 & 0 & 0 & 0 & 0 & 0 & \ldots   \\
   1 & 1 & 0 & 0 & 0 & 0 & 0 & \ldots   \\
   1 & 0 & 1 & 0 & 0 & 0 & 0 & \ldots   \\
   1 & 1 & 1 & 1 & 0 & 0 & 0 & \ldots   \\
   1 & 0 & 2 & 0 & 1 & 0 & 0 & \ldots   \\
   1 & 1 & 2 & 2 & 1 & 1 & 0 & \ldots   \\
   1 & 0 & 3 & 0 & 3 & 0 & 1 & \ldots   \\
   \vdots  & \vdots  & \vdots  & \vdots  & \vdots  & \vdots  & \vdots  & \ddots   \\
\end{matrix} \right), \qquad\left( \begin{matrix}
   1 & 0 & 0 & 0 & 0 & 0 & 0 & \ldots   \\
   -1 & 1 & 0 & 0 & 0 & 0 & 0 & \ldots   \\
   -1 & 0 & 1 & 0 & 0 & 0 & 0 & \ldots   \\
   1 & -1 & -1 & 1 & 0 & 0 & 0 & \ldots   \\
   1 & 0 & -2 & 0 & 1 & 0 & 0 & \ldots   \\
   -1 & 1 & 2 & -2 & -1 & 1 & 0 & \ldots   \\
   -1 & 0 & 3 & 0 & -3 & 0 & 1 & \ldots   \\
   \vdots  & \vdots  & \vdots  & \vdots  & \vdots  & \vdots  & \vdots  & \ddots   \\
\end{matrix} \right),$$
where the series $g\left( -1,x \right)$  is not defined. Since
$${{\left( {{P}_{g\left( -1,x \right)}} \right)}_{2n+i,2m+j}}=\left[ {{x}^{2n+i}} \right]\frac{{{\left( 1+x \right)}^{1-j}}{{x}^{2m+j}}}{{{\left( 1-{{x}^{2}} \right)}^{m+1}}}=\left( \begin{matrix}
   n  \\
   m  \\
\end{matrix} \right), \qquad i\ge j; \qquad=0, \qquad i<j;$$ 
$i,j=0,$$1$, then
$${{P}_{g\left( -1,x \right)}}={}_{0,2}P\times {{P}_{c\left( x \right)}},  \qquad c\left( x \right)=\left( 1+x \right){{e}^{{{x}^{2}}}}:$$
$${{c}_{2n+i}}=\frac{1}{n!},  \qquad 0\le i<2;  \qquad{{c}_{2n-i}}=\frac{1}{\left( n-1 \right)!},  \qquad 0<i\le 2,$$
$${{\left( {{P}_{c\left( x \right)}} \right)}_{2n+i,2m+j}}=\left( \begin{matrix}
   n  \\
   m  \\
\end{matrix} \right), \qquad i\ge j; \qquad=n\left( \begin{matrix}
   n-1  \\
   m  \\
\end{matrix} \right), \qquad i<j.$$

Each nonzero generalized Pascal matrix is the Hadamard product of the matrices $_{\varphi ,q}P$. Since the first column of the matrix ${{P}_{c\left( x \right)}}$, –  denote it $b\left( x \right)$, –  is the Hadamard product of  the first columns of the  matrices $_{\varphi ,q}P$, –  denote them $_{\varphi ,q}b\left( x \right)$:
$$\left[ {{x}^{n}} \right]{}_{\varphi ,q}b\left( x \right)=1, \qquad n\left( \bmod q \right)\ne 0; \qquad=\varphi , \qquad n\left( \bmod q \right)=0,$$
then
$${{P}_{c\left( x \right)}}={}_{2}P\left( {{b}_{2}} \right)\times {}_{3}P\left( {{b}_{3}} \right)\times {}_{4}P\left( {{{b}_{4}}}/{{{b}_{2}}}\; \right)\times {}_{5}P\left( {{b}_{5}} \right)\times {}_{6}P\left( {{{b}_{6}}}/{{{b}_{2}}{{b}_{3}}}\; \right)\times {}_{7}P\left( {{b}_{7}} \right)\times $$
$${{\times }_{8}}P\left( {{{b}_{8}}}/{{{b}_{4}}}\; \right)\times {}_{9}P\left( {{{b}_{9}}}/{{{b}_{3}}}\; \right)\times {}_{10}P\left( {{{b}_{10}}}/{{{b}_{2}}{{b}_{5}}}\; \right)\times {}_{11}P\left( {{b}_{11}} \right)\times {}_{12}P\left( {{{b}_{12}}{{b}_{2}}}/{{{b}_{4}}{{b}_{6}}}\; \right)\times ...$$ 
and so on. Let ${{e}_{q}}$ is a basis vector of an infinite-dimensional vector space. Mapping of the set of generalized Pascal matrices in an infinite-dimensional vector space such that $_{\varphi ,q}P\to {{e}_{q}}\log \left| \varphi  \right|$ is a group  homomorphism whose kernel consists of all involutions in the group of generalized Pascal matrices, i.e. from matrices whose non-zero elements equal to $\pm 1$. Thus, the set of generalized Pascal matrices whose elements are non-negative numbers is an infinite-dimensional vector space. Zero generalized Pascal matrices can be viewed as points at infinity of space.
\section{Algebra associated with zero generalized Pascal matrix}

Let $_{0}P$ is a zero generalized Pascal matrix. Denote
$$\left( a\left( x \right),x \right)\times {}_{0}P=\left( a\left( x \right),x|{}_{0}P \right).$$
Then
$$\left( a\left( x \right),x|{}_{0}P \right)b\left( x \right)=a\left( x \right)\circ b\left( x \right),$$
$$\left( a\left( x \right),x|{}_{0}P \right)\left( b\left( x \right),x|{}_{0}P \right)=\left( a\left( x \right)\circ b\left( x \right),x|{}_{0}P \right),$$
$$\left[ {{x}^{n}} \right]a\left( x \right)\circ b\left( x \right)=\sum\limits_{m=0}^{n}{\left( n,m \right){{a}_{n-m}}{{b}_{m}}}, \qquad\left( n,m \right)={{\left( {}_{0}P \right)}_{n,m}}.$$
Ordinary operation of multiplication of the series remains a priority:
$$a\left( x \right)b\left( x \right)\circ c\left( x \right)=\left( a\left( x \right)b\left( x \right) \right)\circ c\left( x \right).$$
Multiplying the  identity
$$D\left( a\left( x \right),x \right)=\left( a\left( x \right),x \right)D+\left( {a}'\left( x \right),x \right),$$
where$D$ is the matrix of the differential operator, by the matrix $\left( x,x \right)$ we obtain the  identity
$$\left( x,x \right)D\left( a\left( x \right),x \right)=\left( a\left( x \right),x \right)\left( x,x \right)D+\left( x{a}'\left( x \right),x \right),$$
where $\left( x,x \right)D$  is the diagonal matrix. Hence, true the identity
\[\left( x,x \right)D\left( a\left( x \right),x|{}_{0}P \right)=\left( a\left( x \right),x|{}_{0}P \right)\left( x,x \right)D+\left( x{a}'\left( x \right),x|{}_{0}P \right).\]
Thus,
$$x{{\left( a\left( x \right)\circ b\left( x \right) \right)}^{\prime }}=a\left( x \right)\circ x{b}'\left( x \right)+x{a}'\left( x \right)\circ b\left( x \right).$$
Denote
$${{a}^{\left( n \right)}}\left( x \right)=a\left( x \right)\circ {{a}^{\left( n-1 \right)}}\left( x \right),  \qquad{{a}^{\left( 0 \right)}}\left( x \right)=1;$$
then       
$$x{{\left( {{a}^{\left( n \right)}}\left( x \right) \right)}^{\prime }}=n{{a}^{\left( n-1 \right)}}\left( x \right)\circ x{a}'\left( x \right).\eqno(3)$$
Obviously, if
 $$\left( \sum\limits_{n=0}^{\infty }{{{b}_{n}}{{x}^{n}}} \right)\left( \sum\limits_{n=0}^{\infty }{{{c}_{n}}{{x}^{n}}} \right)=\sum\limits_{n=0}^{\infty }{{{d}_{n}}{{x}^{n}}},$$
then
$$\left( \sum\limits_{n=0}^{\infty }{{{b}_{n}}{{a}^{\left( n \right)}}\left( x \right)} \right)\circ \left( \sum\limits_{n=0}^{\infty }{{{c}_{n}}{{a}^{\left( n \right)}}\left( x \right)} \right)=\sum\limits_{n=0}^{\infty }{{{d}_{n}}}{{a}^{\left( n \right)}}\left( x \right).$$
Therefore the power and the logarithm of series  defined the same way as in the ordinary algebra of formal power series:
$${{a}^{\left( \varphi  \right)}}\left( x \right)=\sum\limits_{n=0}^{\infty }{\left( \begin{matrix}
   \varphi   \\
   n  \\
\end{matrix} \right){{\left( a\left( x \right)-1 \right)}^{\left( n \right)}}}, \qquad{{a}_{0}}=1,$$
$$\log \circ a\left( x \right)=\sum\limits_{n=1}^{\infty }{\frac{{{\left( -1 \right)}^{n-1}}}{n}}{{\left( a\left( x \right)-1 \right)}^{\left( n \right)}},$$
$${{a}^{\left( \varphi  \right)}}\left( x \right)=\sum\limits_{n=0}^{\infty }{\frac{{{\varphi }^{n}}}{n!}}{{\left( \log \circ a\left( x \right) \right)}^{\left( n \right)}}.$$
Then
$${{a}^{\left( \varphi  \right)}}\left( x \right)\circ {{a}^{\left( \beta  \right)}}\left( x \right)={{a}^{\left( \varphi +\beta  \right)}}\left( x \right), \qquad{{a}^{\left( \varphi  \right)}}\left( x \right)\circ {{b}^{\left( \varphi  \right)}}\left( x \right)={{\left( a\left( x \right)\circ b\left( x \right) \right)}^{\left( \varphi  \right)}},$$
$$\log \circ {{a}^{\left( \varphi  \right)}}\left( x \right)=\varphi \log \circ a\left( x \right),\qquad\log \circ \left( a\left( x \right)\circ b\left( x \right) \right)=\log \circ a\left( x \right)+\log \circ b\left( x \right).$$
Note the identity
	$$x{{\left( \log \circ a\left( x \right) \right)}^{\prime }}=x{a}'\left( x \right)\circ {{a}^{\left( -1 \right)}}\left( x \right).	\eqno(4)$$

\section{Basic zero generalized Pascal matrices }

Basic zero generalized Pascal matrices, as is clear from the definition, are the matrices $_{0,q}P$,
$${{\left( {}_{0,q}P \right)}_{n,m}}=1, \qquad n\left( \bmod q \right)\ge m\left( \bmod q \right); \qquad=0, \qquad n\left( \bmod q \right)<m\left( \bmod q \right).$$
For example, $_{0,2}P$, $_{0,3}P$:
$$\left( \begin{matrix}
   1 & 0 & 0 & 0 & 0 & 0 & 0 & 0 & 0 & \ldots   \\
   1 & 1 & 0 & 0 & 0 & 0 & 0 & 0 & 0 & \ldots   \\
   1 & 0  & 1 & 0 & 0 & 0 & 0 & 0 & 0 & \ldots   \\
   1 & 1 & 1 & 1 & 0 & 0 & 0 & 0 & 0 & \ldots   \\
   1 & 0  & 1 & 0  & 1 & 0 & 0 & 0 & 0 & \ldots   \\
   1 & 1 & 1 & 1 & 1 & 1 & 0 & 0 & 0 & \dots   \\
   1 & 0  & 1 & 0  & 1 & 0  & 1 & 0 & 0 & \ldots   \\
   1 & 1 & 1 & 1 & 1 & 1 & 1 & 1 & 0 & \ldots   \\
   1 & 0  & 1 & 0  & 1 & 0 & 1 & 0  & 1 & \ldots   \\
   \vdots  & \vdots  & \vdots  & \vdots  & \vdots  & \vdots  & \vdots  & \vdots  & \vdots  & \ddots   \\
\end{matrix} \right), \qquad \left( \begin{matrix}
   1 & 0 & 0 & 0 & 0 & 0 & 0 & 0 & 0 & \ldots   \\
   1 & 1 & 0 & 0 & 0 & 0 & 0 & 0 & 0 & \ldots   \\
   1 & 1 & 1 & 0 & 0 & 0 & 0 & 0 & 0 & \ldots   \\
   1 & 0  & 0 & 1 & 0 & 0 & 0 & 0 & 0 & \ldots   \\
   1 & 1 & 0  & 1 & 1 & 0 & 0 & 0 & 0 & \ldots   \\
   1 & 1 & 1 & 1 & 1 & 1 & 0 & 0 & 0 & \ldots   \\
   1 & 0  & 0  & 1 & 0  & 0 & 1 & 0 & 0 & \ldots   \\
   1 & 1 & 0  & 1 & 1 & 0  & 1 & 1 & 0 & \ldots   \\
   1 & 1 & 1 & 1 & 1 & 1 & 1 & 1 & 1 & \ldots   \\
   \vdots  & \vdots  & \vdots  & \vdots  & \vdots  & \vdots  & \vdots  & \vdots  & \vdots  & \ddots   \\
\end{matrix} \right).$$
Note some characteristic properties of these matrices and associated algebras.

The block matrix, whose $\left( n,m \right)$th block  is the matrix consisting of $q$ first rows of the matrix $\left( b\left( x \right),x \right)$ and multiplied by ${{a}_{n-m}}$, denote ${{\left( a\left( x \right)|b\left( x \right) \right)}_{q}}$. For example,
$${{\left( a\left( x \right)|b\left( x \right) \right)}_{2}}=\left( \begin{matrix}
   {{a}_{0}}{{b}_{0}} & 0 & 0 & 0 & 0 & 0 & \ldots   \\
   {{a}_{0}}{{b}_{1}} & {{a}_{0}}{{b}_{0}} & 0 & 0 & 0 & 0 & \ldots   \\
   {{a}_{1}}{{b}_{0}} & 0 & {{a}_{0}}{{b}_{0}} & 0 & 0 & 0 & \ldots   \\
   {{a}_{1}}{{b}_{1}} & {{a}_{1}}{{b}_{0}} & {{a}_{0}}{{b}_{1}} & {{a}_{0}}{{b}_{0}} & 0 & 0 & \ldots   \\
   {{a}_{2}}{{b}_{0}} & 0 & {{a}_{1}}{{b}_{0}} & 0 & {{a}_{0}}{{b}_{0}} & 0 & \ldots   \\
   {{a}_{2}}{{b}_{1}} & {{a}_{2}}{{b}_{0}} & {{a}_{1}}{{b}_{1}} & {{a}_{1}}{{b}_{0}} & {{a}_{0}}{{b}_{1}} & {{a}_{0}}{{b}_{0}} & \ldots   \\
   \vdots  & \vdots  & \vdots  & \vdots  & \vdots  & \vdots  & \ddots   \\
\end{matrix} \right),$$
$${{\left( a\left( x \right)|b\left( x \right) \right)}_{3}}=\left( \begin{matrix}
   {{a}_{0}}{{b}_{0}} & 0 & 0 & 0 & 0 & 0 & \ldots   \\
   {{a}_{0}}{{b}_{1}} & {{a}_{0}}{{b}_{0}} & 0 & 0 & 0 & 0 & \ldots   \\
   {{a}_{0}}{{b}_{2}} & {{a}_{0}}{{b}_{1}} & {{a}_{0}}{{b}_{0}} & 0 & 0 & 0 & \ldots   \\
   {{a}_{1}}{{b}_{0}} & 0 & 0 & {{a}_{0}}{{b}_{0}} & 0 & 0 & \ldots   \\
   {{a}_{1}}{{b}_{1}} & {{a}_{1}}{{b}_{0}} & 0 & {{a}_{0}}{{b}_{1}} & {{a}_{0}}{{b}_{0}} & 0 & \ldots   \\
   {{a}_{1}}{{b}_{2}} & {{a}_{1}}{{b}_{1}} & {{a}_{1}}{{b}_{0}} & {{a}_{0}}{{b}_{2}} & {{a}_{0}}{{b}_{1}} & {{a}_{0}}{{b}_{0}} & \ldots   \\
   \vdots  & \vdots  & \vdots  & \vdots  & \vdots  & \vdots  & \ddots   \\
\end{matrix} \right).$$
Then
$${{\left( a\left( x \right)|b\left( x \right) \right)}_{q}}=\left( \left( \sum\limits_{n=0}^{q-1}{{{b}_{n}}{{x}^{n}}} \right)a\left( {{x}^{q}} \right),x|{}_{0,q}P \right),$$
$${{\left( a\left( x \right)|b\left( x \right) \right)}_{q}}{{\left( g\left( x \right)|c\left( x \right) \right)}_{q}}={{\left( a\left( x \right)g\left( x \right)|b\left( x \right)c\left( x \right) \right)}_{q}},$$
$$_{0,q}{{P}^{\varphi }}={{\left( {{\left( \frac{1}{1-x} \right)}^{\varphi }}|{{\left( \frac{1}{1-x} \right)}^{\varphi }} \right)}_{q}}=\left( {{\left( \frac{1}{1-x} \right)}^{\left( \varphi  \right)}},x|{}_{0,q}P \right),$$
$${{\left( \frac{1}{1-x} \right)}^{\left( \varphi  \right)}}=\left( \sum\limits_{m=0}^{q-1}{\left( \begin{matrix}
   \varphi +m-1  \\
   m  \\
\end{matrix} \right){{x}^{m}}} \right){{\left( \frac{1}{1-{{x}^{q}}} \right)}^{\varphi }},$$
$$\left[ {{x}^{qn+i}} \right]{{\left( \frac{1}{1-x} \right)}^{\left( \varphi  \right)}}=\left( \begin{matrix}
   \varphi +n-1  \\
   n  \\
\end{matrix} \right)\left( \begin{matrix}
   \varphi +i-1  \\
   i  \\
\end{matrix} \right), \qquad0\le i<q,$$
$${{\left( _{0,q}{{P}^{\varphi }} \right)}_{qn+i,qm+j}}=\left( \begin{matrix}
   \varphi +n-m-1  \\
   n-m  \\
\end{matrix} \right)\left( \begin{matrix}
   \varphi +i-j-1  \\
   i-j  \\
\end{matrix} \right); \qquad\left( \begin{matrix}
   \beta   \\
   k  \\
\end{matrix} \right)=0, \qquad k<0.$$

Denote
\[{{w}_{n}}\left( \varphi ,x \right)=\sum\limits_{m=0}^{n}{\left( \begin{matrix}
   \varphi +m-1  \\
   m  \\
\end{matrix} \right){{x}^{n-m}}}.\]
Then
$$\left[ qn+m,\to  \right]{}_{0,q}{{P}^{\varphi }}={{w}_{m}}\left( \varphi ,x \right){{w}_{n}}\left( \varphi ,{{x}^{q}} \right),  \qquad0\le m<q.$$
Based on the formula (4) we find the series $\log \circ {{\left( 1-x \right)}^{-1}}$:
$$x{{\left( \log \circ {{\left( 1-x \right)}^{-1}} \right)}^{\prime }}=x{{\left( 1-x \right)}^{-2}}\circ {{\left( \frac{1}{1-x} \right)}^{\left( -1 \right)}}={}_{0,q}{{P}^{-1}}x{{\left( 1-x \right)}^{-2}}.$$
Since
$${{w}_{0}}\left( -1,x \right)=1,  \qquad{{w}_{n}}\left( -1,x \right)={{x}^{n-1}}\left( x-1 \right),$$
then
$$\left[ qn+m,\to  \right]{}_{0,q}{{P}^{-1}}={{x}^{m-1}}\left( x-1 \right), \qquad n=0, \qquad m\ne 0;$$ 
$$={{x}^{q\left( n-1 \right)}}\left( {{x}^{q}}-1 \right), \quad m=0, \quad n\ne 0; =\quad{{x}^{q\left( n-1 \right)+m-1}}\left( x-1 \right)\left( {{x}^{q}}-1 \right), \quad n,m\ne 0.$$
For example,
$$_{0,3}{{P}^{-1}}=\left( \begin{matrix}
   1 & 0 & 0 & 0 & 0 & 0 & 0 & 0 & 0 & \ldots   \\
   -1 & 1 & 0 & 0 & 0 & 0 & 0 & 0 & 0 & \ldots   \\
   0 & -1 & 1 & 0 & 0 & 0 & 0 & 0 & 0 & \ldots   \\
   -1 & 0 & 0 & 1 & 0 & 0 & 0 & 0 & 0 & \ldots   \\
   1 & -1 & 0 & -1 & 1 & 0 & 0 & 0 & 0 & \ldots   \\
   0 & 1 & -1 & 0 & -1 & 1 & 0 & 0 & 0 & \ldots   \\
   0 & 0 & 0 & -1 & 0 & 0 & 1 & 0 & 0 & \ldots   \\
   0 & 0 & 0 & 1 & -1 & 0 & -1 & 1 & 0 & \ldots   \\
   0 & 0 & 0 & 0 & 1 & -1 & 0 & -1 & \text{ }1 & \ldots   \\
   \vdots  & \vdots  & \vdots  & \vdots  & \vdots  & \vdots  & \vdots  & \vdots  & \vdots  & \ddots   \\
\end{matrix} \right).$$
We use the theorem that will be needed in the more complex cases: if a polynomial has the form
$$\prod\limits_{i=1}^{n}{\left( {{x}^{{{m}_{i}}}}-1 \right)}={{\left( -1 \right)}^{n}}\left( 1+\sum\limits_{i=1}^{{{2}^{n-1}}-1}{{{x}^{{{p}_{i}}}}}-\sum\limits_{i=1}^{{{2}^{n-1}}}{{{x}^{{{s}_{i}}}}} \right),  \qquad n>1,$$
then
$$\sum\limits_{i=1}^{{{2}^{n-1}}-1}{{{p}_{i}}}=\sum\limits_{i=1}^{{{2}^{n-1}}}{{{s}_{i}}}.$$
This follows from the Vieta’s formulas, if we consider the monomials ${{x}^{{{m}_{i}}}}$ as the roots of a polynomial:
$$\prod\limits_{i=1}^{n}{\left( 1-{{x}^{{{m}_{i}}}} \right)}=1-\left( {{x}^{{{m}_{1}}}}+{{x}^{{{m}_{2}}}}+...+{{x}^{{{m}_{n}}}} \right)+\left( {{x}^{{{m}_{1}}+{{m}_{2}}}}+{{x}^{{{m}_{1}}+{{m}_{3}}}}+...+{{x}^{{{m}_{n-1}}+{{m}_{n}}}} \right)-$$ 
$$-\left( {{x}^{{{m}_{1}}+{{m}_{2}}+{{m}_{3}}}}+{{x}^{{{m}_{1}}+{{m}_{2}}+{{m}_{4}}}}+...+{{x}^{{{m}_{n-2}}+{{m}_{n-1}}+{{m}_{n}}}} \right)+...{{\left( -1 \right)}^{n}}{{x}^{{{m}_{1}}+{{m}_{2}}...+{{m}_{n}}}}.$$
Hence, if  the $n$th row of the matrix $A$ has the form
${{x}^{p}}\prod\nolimits_{i=1}^{s}{\left( {{x}^{{{m}_{i}}}}-1 \right)}$, $s>1$, then $\left[ {{x}^{n}} \right]Ax{{\left( 1-x \right)}^{-2}}=0$. Thus,
$$_{0,q}{{P}^{-1}}x{{\left( 1-x \right)}^{-2}}=\sum\limits_{m=1}^{q-1}{{{x}^{m}}}+\sum\limits_{m=1}^{\infty }{q{{x}^{qm}}},$$
$$\log \circ {{\left( 1-x \right)}^{-1}}=\sum\limits_{m=1}^{q-1}{\frac{{{x}^{m}}}{m}}+\log {{\left( 1-{{x}^{q}} \right)}^{-1}}.$$

Attribute of the algebras associated with zero generalized Pascal matrices are the series of the form $a\left( x \right)=1+\log \circ a\left( x \right)$ that satisfy the identities
  $$\left( \log \circ a\left( x \right) \right)\circ \left( \log \circ a\left( x \right) \right)=0,$$
$${{a}^{\left( \varphi  \right)}}\left( x \right)=\sum\limits_{n=0}^{\infty }{\frac{{{\varphi }^{n}}}{n!}}{{\left( \log \circ a\left( x \right) \right)}^{\left( n \right)}}=1+\varphi \log \circ a\left( x \right).$$
We call these series, for example, $l$-series. In the algebra associated with the matrix $_{0,q}P$, $l$-series form the group whose elements are multiplied by the rule
$${{a}_{1}}\left( x \right)\circ {{a}_{2}}\left( x \right)=1+\log \circ \left( {{a}_{1}}\left( x \right)\circ {{a}_{2}}\left( x \right) \right).$$
We will find a general view of them as follows. As seen from the matrix $_{0,2}P$, monomials ${{x}^{2n+1}}$ form a closed system of zero divisors, i. e. their products with each other and with ourselves are zero. Therefore, $l$-series has the form $a\left( x \right)=1+xb\left( {{x}^{2}} \right)$. As seen from the matrix $_{0,3}P$, monomials ${{x}^{3n+2}}$ form a closed system of zero divisors. Therefore, $l$-series has the form $a\left( x \right)=1+{{x}^{2}}b\left( {{x}^{3}} \right)$. As seen from the matrix 
$$_{0,4}P=\left( \begin{matrix}
   1 & 0 & 0 & 0 & 0 & 0 & 0 & 0 & \ldots   \\
   1 & 1 & 0 & 0 & 0 & 0 & 0 & 0 & \ldots   \\
   1 & 1 & 1 & 0 & 0 & 0 & 0 & 0 & \ldots   \\
   1 & 1 & 1 & 1 & 0 & 0 & 0 & 0 & \ldots   \\
   1 & 0 & 0 & 0 & 1 & 0 & 0 & 0 & \ldots   \\
   1 & 1 & 0 & 0 & 1 & 1 & 0 & 0 & \ldots   \\
   1 & 1 & 1 & 0 & 1 & 1 & 1 & 0 & \ldots   \\
   1 & 1 & 1 & 1 & 1 & 1 & 1 & 1 & \ldots   \\
   \vdots  & \vdots  & \vdots  & \vdots  & \vdots  & \vdots  & \vdots  & \vdots  & \ddots   \\
\end{matrix} \right),$$
monomials ${{x}^{4n+2}}$, ${{x}^{4n+3}}$ form a closed system of zero divisors. Therefore, $l$-series has the form $a\left( x \right)=1+{{x}^{2}}{{b}_{1}}\left( {{x}^{4}} \right)+{{x}^{3}}{{b}_{2}}\left( {{x}^{4}} \right)$. In general, in the algebra associated with the matrix $_{0,q}P$, monomials ${{x}^{qn+q-m}}$, $1\le m\le \left\lfloor {q}/{2}\; \right\rfloor $, where $\left\lfloor {q}/{2}\; \right\rfloor $ is the integral part of ${q}/{2}\;$, form a closed system of zero divisors;  $l$-series has the form 
$$a\left( x \right)=1+\sum\limits_{m=1}^{\left\lfloor {q}/{2}\; \right\rfloor }{{{x}^{q-m}}}{{b}_{\left\lfloor {q}/{2}\; \right\rfloor -m+1}}\left( {{x}^{q}} \right).\eqno(5)$$
Evident that the algebra associated with the Hadamard product of the  matrices $_{0,q}P$ contains the all  groups of $l$-series of the algebras associated with the factors. For example, in the algebra associated with the matrix
$${}_{0,2}P\times {}_{0,3}P\times {}_{0,4}P\times \text{ }...\text{ =}\left( \begin{matrix}
   1 & 0 & 0 & 0 & 0 & \ldots   \\
   1 & 1 & 0 & 0 & 0 & \ldots   \\
   1 & 0 & 1 & 0 & 0 & \ldots   \\
   1 & 0 & 0 & 1 & 0 & \ldots   \\
   1 & 0 & 0 & 0 & 1 & \ldots   \\
   \vdots  & \vdots  & \vdots  & \vdots  & \vdots  & \ddots   \\
\end{matrix} \right)$$
$l$-series are the all series with ${{a}_{0}}=1$.

\section{Exponential zero generalized Pascal  matrices}

Consider matrix
$${{P}_{c\left( q,x \right)}}, \qquad c\left( q,x \right)=\left( \sum\limits_{n=0}^{q-1}{{{x}^{n}}} \right){{e}^{{{x}^{q}}}},$$
$${{c}_{qn+i}}=\frac{1}{n!},  \qquad0\le i<q;  \qquad{{c}_{qn-i}}=\frac{1}{\left( n-1 \right)!},  \qquad0<i\le q,$$
$${{\left( {{P}_{c\left( q,x \right)}} \right)}_{qn+i,qm+j}}=\left( \begin{matrix}
   n  \\
   m  \\
\end{matrix} \right), \qquad i\ge j; \qquad=n\left( \begin{matrix}
   n-1  \\
   m  \\
\end{matrix} \right), i<j.$$
Matrix ${}_{0,q}{{P}_{e}}={}_{0,q}P\times {{P}_{c\left( q,x \right)}}$ will be called exponential zero generalized Pascal  matrix. For example, ${}_{0,2}{{P}_{e}}={{P}_{g\left( -1,x \right)}}$; matrices  ${}_{0,3}{{P}_{e}}$, ${}_{0,3}P_{e}^{-1}$ have the form
$$\left( \begin{matrix}
   1 & 0 & 0 & 0 & 0 & 0 & 0 & 0 & 0 & \ldots   \\
   1 & 1 & 0 & 0 & 0 & 0 & 0 & 0 & 0 & \ldots   \\
   1 & 1 & 1 & 0 & 0 & 0 & 0 & 0 & 0 & \ldots   \\
   1 & 0 & 0 & 1 & 0 & 0 & 0 & 0 & 0 & \ldots   \\
   1 & 1 & 0 & 1 & 1 & 0 & 0 & 0 & 0 & \ldots   \\
   1 & 1 & 1 & 1 & 1 & 1 & 0 & 0 & 0 & \ldots   \\
   1 & 0 & 0 & 2 & 0 & 0 & 1 & 0 & 0 & \ldots   \\
   1 & 1 & 0 & 2 & 2 & 0 & 1 & 1 & 0 & \ldots   \\
   1 & 1 & 1 & 2 & 2 & 2 & 1 & 1 & 1 & \ldots   \\
   \vdots  & \vdots  & \vdots  & \vdots  & \vdots  & \vdots  & \vdots  & \vdots  & \vdots  & \ddots   \\
\end{matrix} \right),  \quad\left( \begin{matrix}
   1 & 0 & 0 & 0 & 0 & 0 & 0 & 0 & 0 & \ldots   \\
   -1 & 1 & 0 & 0 & 0 & 0 & 0 & 0 & 0 & \ldots   \\
   0 & -1 & 1 & 0 & 0 & 0 & 0 & 0 & 0 & \ldots   \\
   -1 & 0 & 0 & 1 & 0 & 0 & 0 & 0 & 0 & \ldots   \\
   1 & -1 & 0 & -1 & 1 & 0 & 0 & 0 & 0 & \ldots   \\
   0 & 1 & -1 & 0 & -1 & 1 & 0 & 0 & 0 & \ldots   \\
   1 & 0 & 0 & -2 & 0 & 0 & 1 & 0 & 0 & \ldots   \\
   -1 & 1 & 0 & 2 & -2 & 0 & -1 & 1 & 0 & \ldots   \\
   0 & -1 & 1 & 0 & 2 & -2 & 0 & -1 & 1 & \ldots   \\
   \vdots  & \vdots  & \vdots  & \vdots  & \vdots  & \vdots  & \vdots  & \vdots  & \vdots  & \ddots   \\
\end{matrix} \right).$$

The block matrix, whose $\left( n,m \right)$th block  is the matrix consisting of $q$ first rows of the matrix $\left( b\left( x \right),x \right)$ and multiplied by
${{a}_{n-m}}{n\choose m}$, denote ${{\left( a\left( x \right)|b\left( x \right) \right)}_{q,e}}$ . Then
$${{\left( a\left( x \right)|b\left( x \right) \right)}_{q,e}}=\left( \left( \sum\limits_{n=0}^{q-1}{{{b}_{n}}{{x}^{n}}} \right)a\left( {{x}^{q}} \right),x|{}_{0,q}{{P}_{e}} \right),$$
$${{\left( a\left( x \right)|b\left( x \right) \right)}_{q,e}}{{\left( g\left( x \right)|c\left( x \right) \right)}_{q,e}}={{\left( a\left( x \right)\circ g\left( x \right)|b\left( x \right)c\left( x \right) \right)}_{q,e}},$$
where
$$\left[ {{x}^{n}} \right]a\left( x \right)\circ g\left( x \right)=\sum\limits_{m=0}^{n}{\left( \begin{matrix}
   n  \\
   m  \\
\end{matrix} \right){{a}_{m}}{{g}_{n-m}}},$$
$$_{0,q}P_{e}^{\varphi }={{\left( \left( \frac{1}{1-\varphi x} \right)|{{\left( \frac{1}{1-x} \right)}^{\varphi }} \right)}_{q,e}}=\left( {{\left( \frac{1}{1-x} \right)}^{\left( \varphi  \right)}},x|{}_{0,q}{{P}_{e}} \right),$$
$${{\left( \frac{1}{1-x} \right)}^{\left( \varphi  \right)}}=\left( \sum\limits_{m=0}^{q-1}{\left( \begin{matrix}
   \varphi +m-1  \\
   m  \\
\end{matrix} \right){{x}^{m}}} \right)\frac{1}{1-\varphi {{x}^{q}}},$$
$$\left[ {{x}^{qn+i}} \right]{{\left( \frac{1}{1-x} \right)}^{\left( \varphi  \right)}}={{\varphi }^{n}}\left( \begin{matrix}
   \varphi +i-1  \\
   i  \\
\end{matrix} \right), \qquad0\le i<q,$$
$${{\left( _{0,q}P_{e}^{\varphi } \right)}_{qn+i,qm+j}}={{\varphi }^{n-m}}\left( \begin{matrix}
   \varphi +i-j-1  \\
   i-j  \\
\end{matrix} \right)\left( \begin{matrix}
   n  \\
   m  \\
\end{matrix} \right),$$ 
$$\left[ qn+m,\to  \right]{}_{0,q}P_{e}^{\varphi }={{w}_{m}}\left( \varphi ,x \right){{\left( {{x}^{q}}+\varphi  \right)}^{n}},  \qquad0\le m<q.$$
Since
$$\left[ qn+m,\to  \right]{}_{0,q}P_{e}^{-1}={{x}^{m-1}}\left( x-1 \right){{\left( {{x}^{q}}-1 \right)}^{n}}, \quad m>0; \quad={{\left( {{x}^{q}}-1 \right)}^{n}}, \quad m=0,$$
then
$$_{0,q}{{P}^{-1}}x{{\left( 1-x \right)}^{-2}}=\sum\limits_{m=1}^{q-1}{{{x}^{m}}}+q{{x}^{q}},$$
$$\log \circ {{\left( 1-x \right)}^{-1}}=\sum\limits_{m=1}^{q-1}{\frac{{{x}^{m}}}{m}}+\log \left( \exp {{x}^{q}} \right).$$

\section{Fractal zero generalized Pascal matrices}

Consider matrix

\[_{\left[ 0,q \right]}P={}_{0,q}P\times {}_{0,{{q}^{2}}}P\times {}_{0,{{q}^{3}}}P\times ...\times {}_{0,{{q}^{k}}}P\times ...\] 
For example (Pascal triangle modulo 2),

$$_{\setcounter{MaxMatrixCols}{20}\left[ 0,2 \right]}P=\left( \begin{matrix}
   1 & 0 & 0 & 0 & 0 & 0 & 0 & 0 & 0 & 0 & 0 & 0 & 0 & 0 & 0 & 0 & \ldots   \\
   1 & 1 & 0 & 0 & 0 & 0 & 0 & 0 & 0 & 0 & 0 & 0 & 0 & 0 & 0 & 0 & \ldots   \\
   1 & 0 & 1 & 0 & 0 & 0 & 0 & 0 & 0 & 0 & 0 & 0 & 0 & 0 & 0 & 0 & \ldots   \\
   1 & 1 & 1 & 1 & 0 & 0 & 0 & 0 & 0 & 0 & 0 & 0 & 0 & 0 & 0 & 0 & \ldots   \\
   1 & 0 & 0 & 0 & 1 & 0 & 0 & 0 & 0 & 0 & 0 & 0 & 0 & 0 & 0 & 0 & \ldots   \\
   1 & 1 & 0 & 0 & 1 & 1 & 0 & 0 & 0 & 0 & 0 & 0 & 0 & 0 & 0 & 0 & \ldots   \\
   1 & 0 & 1 & 0 & 1 & 0 & 1 & 0 & 0 & 0 & 0 & 0 & 0 & 0 & 0 & 0 & \ldots   \\
   1 & 1 & 1 & 1 & 1 & 1 & 1 & 1 & 0 & 0 & 0 & 0 & 0 & 0 & 0 & 0 & \ldots   \\
   1 & 0 & 0 & 0 & 0 & 0 & 0 & 0 & 1 & 0 & 0 & 0 & 0 & 0 & 0 & 0 & \ldots   \\
   1 & 1 & 0 & 0 & 0 & 0 & 0 & 0 & 1 & 1 & 0 & 0 & 0 & 0 & 0 & 0 & \ldots   \\
   1 & 0 & 1 & 0 & 0 & 0 & 0 & 0 & 1 & 0 & 1 & 0 & 0 & 0 & 0 & 0 & \ldots   \\
   1 & 1 & 1 & 1 & 0 & 0 & 0 & 0 & 1 & 1 & 1 & 1 & 0 & 0 & 0 & 0 & \ldots   \\
   1 & 0 & 0 & 0 & 1 & 0 & 0 & 0 & 1 & 0 & 0 & 0 & 1 & 0 & 0 & 0 & \ldots   \\
   1 & 1 & 0 & 0 & 1 & 1 & 0 & 0 & 1 & 1 & 0 & 0 & 1 & 1 & 0 & 0 & \ldots   \\
   1 & 0 & 1 & 0 & 1 & 0 & 1 & 0 & 1 & 0 & 1 & 0 & 1 & 0 & 1 & 0 & \ldots   \\
   1 & 1 & 1 & 1 & 1 & 1 & 1 & 1 & 1 & 1 & 1 & 1 & 1 & 1 & 1 & 1 & \ldots   \\
   \vdots  & \vdots  & \vdots  & \vdots  & \vdots  & \vdots  & \vdots  & \vdots  & \vdots  & \vdots  & \vdots  & \vdots  & \vdots  & \vdots  & \vdots  & \vdots  & \ddots   \\
\end{matrix} \right).$$
$${{\left( {}_{\left[ 0,q \right]}P \right)}_{n,m}}=1, \quad n\left( \bmod {{q}^{k}} \right)\ge m\left( \bmod {{q}^{k}} \right);\quad=0, \quad n\left( \bmod {{q}^{k}} \right)<m\left( \bmod {{q}^{k}} \right), \quad k=1, 2, … .$$
Denote
$${{\left( _{\left[ 0,q \right]}P \right)}_{n,m}}={{\left( \begin{matrix}
   n  \\
   m  \\
\end{matrix} \right)}_{0,q}}.$$
{\bfseries Theorem.}
$${{\left( \begin{matrix}
   {{q}^{k}}n+i  \\
   {{q}^{k}}m+j  \\
\end{matrix} \right)}_{0,q}}={{\left( \begin{matrix}
   n  \\
   m  \\
\end{matrix} \right)}_{0,q}}{{\left( \begin{matrix}
   i  \\
   j  \\
\end{matrix} \right)}_{0,q}},  \qquad 0\le i,j<{{q}^{k}}.$$ 
{\bfseries Proof.} By definition, if $n\left( \bmod {{q}^{k}} \right)<m\left( \bmod {{q}^{k}} \right)$ for some value of $k$, then ${{n\choose m}_{0,q}}=0$. We represent the numbers $n$, $m$ in the form 
$$n=\sum\limits_{i=0}^{\infty }{{{n}_{i}}}{{q}^{i}},  \qquad m=\sum\limits_{i=0}^{\infty }{{{m}_{i}}{{q}^{i}}}, \qquad 0\le {{n}_{i}},{{m}_{i}}<q.$$
Then
$$n\left( \bmod {{q}^{k}} \right)=\sum\limits_{i=0}^{k-1}{{{n}_{i}}}{{q}^{i}},  \qquad m\left( \bmod {{q}^{k}} \right)=\sum\limits_{i=0}^{k-1}{{{m}_{i}}{{q}^{i}}}.$$
If $n\left( \bmod {{q}^{k}} \right)<m\left( \bmod {{q}^{k}} \right)$, then ${{n}_{i}}<{{m}_{i}}$ at least for one $i$. Since ${{{{n}_{i}}\choose {{m}_{i}}}_{0,q}}=1$, if ${{n}_{i}}\ge {{m}_{i}}$, then true the identity
$${{\left( \begin{matrix}
   n  \\
   m  \\
\end{matrix} \right)}_{0,q}}=\prod\limits_{i=0}^{\infty }{{{\left( \begin{matrix}
   {{n}_{i}}  \\
   {{m}_{i}}  \\
\end{matrix} \right)}_{0,q}}}.$$
It remains to note that if
$$n=\sum\limits_{i=0}^{\infty }{{{n}_{i}}{{q}^{i}}=}{{q}^{k}}s+j, \qquad 0\le j<{{q}^{k}},$$ 
then
$$j=\sum\limits_{i=0}^{k-1}{{{n}_{i}}}{{q}^{i}}, \qquad s=\sum\limits_{i=0}^{\infty }{{{n}_{i+k}}{{q}^{i}}}.$$

Denote
$$\left( \left( \sum\limits_{n=0}^{{{q}^{k}}-1}{{{b}_{n}}{{x}^{n}}} \right)a\left( {{x}^{{{q}^{k}}}} \right),x|{}_{\left[ 0,q \right]}P \right)={{\left( a\left( x \right)|b\left( x \right) \right)}_{q,k}},$$
For exampl,
$${{\left( a\left( x \right)|b\left( x \right) \right)}_{2,1}}=\left( \begin{matrix}
   {{a}_{0}}{{b}_{0}} & 0 & 0 & 0 & 0 & 0 & 0 & 0 & \ldots   \\
   {{a}_{0}}{{b}_{1}} & {{a}_{0}}{{b}_{0}} & 0 & 0 & 0 & 0 & 0 & 0 & \ldots   \\
   {{a}_{1}}{{b}_{0}} & 0 & {{a}_{0}}{{b}_{0}} & 0 & 0 & 0 & 0 & 0 & \ldots   \\
   {{a}_{1}}{{b}_{1}} & {{a}_{1}}{{b}_{0}} & {{a}_{0}}{{b}_{1}} & {{a}_{0}}{{b}_{0}} & 0 & 0 & 0 & 0 & \ldots   \\
   {{a}_{2}}{{b}_{0}} & 0 & 0 & 0 & {{a}_{0}}{{b}_{0}} & 0 & 0 & 0 & \ldots   \\
   {{a}_{2}}{{b}_{1}} & {{a}_{2}}{{b}_{0}} & 0 & 0 & {{a}_{0}}{{b}_{1}} & {{a}_{0}}{{b}_{0}} & 0 & 0 & \ldots   \\
   {{a}_{3}}{{b}_{0}} & 0 & {{a}_{2}}{{b}_{0}} & 0 & {{a}_{1}}{{b}_{0}} & 0 & {{a}_{0}}{{b}_{0}} & 0 & \ldots   \\
   {{a}_{3}}{{b}_{1}} & {{a}_{3}}{{b}_{0}} & {{a}_{2}}{{b}_{1}} & {{a}_{2}}{{b}_{0}} & {{a}_{1}}{{b}_{1}} & {{a}_{1}}{{b}_{0}} & {{a}_{0}}{{b}_{1}} & {{a}_{0}}{{b}_{0}} & \ldots   \\
   \vdots  & \vdots  & \vdots  & \vdots  & \vdots  & \vdots  & \vdots  & \vdots  & \ddots   \\
\end{matrix} \right),$$
$${{\left( a\left( x \right)|b\left( x \right) \right)}_{2,2}}=\left( \begin{matrix}
   {{a}_{0}}{{b}_{0}} & 0 & 0 & 0 & 0 & 0 & 0 & 0 & \ldots   \\
   {{a}_{0}}{{b}_{1}} & {{a}_{0}}{{b}_{0}} & 0 & 0 & 0 & 0 & 0 & 0 & \ldots   \\
   {{a}_{0}}{{b}_{2}} & 0 & {{a}_{0}}{{b}_{0}} & 0 & 0 & 0 & 0 & 0 & \ldots   \\
   {{a}_{0}}{{b}_{3}} & {{a}_{0}}{{b}_{2}} & {{a}_{0}}{{b}_{1}} & {{a}_{0}}{{b}_{0}} & 0 & 0 & 0 & 0 & \ldots   \\
   {{a}_{1}}{{b}_{0}} & 0 & 0 & 0 & {{a}_{0}}{{b}_{0}} & 0 & 0 & 0 & \ldots   \\
   {{a}_{1}}{{b}_{1}} & {{a}_{1}}{{b}_{0}} & 0 & 0 & {{a}_{0}}{{b}_{1}} & {{a}_{0}}{{b}_{0}} & 0 & 0 & \ldots   \\
   {{a}_{1}}{{b}_{2}} & 0 & {{a}_{1}}{{b}_{0}} & 0 & {{a}_{0}}{{b}_{2}} & 0 & {{a}_{0}}{{b}_{0}} & 0 & \ldots   \\
   {{a}_{1}}{{b}_{3}} & {{a}_{1}}{{b}_{2}} & {{a}_{1}}{{b}_{1}} & {{a}_{1}}{{b}_{0}} & {{a}_{0}}{{b}_{3}} & {{a}_{0}}{{b}_{2}} & {{a}_{0}}{{b}_{1}} & {{a}_{0}}{{b}_{0}} & \ldots   \\
   \vdots  & \vdots  & \vdots  & \vdots  & \vdots  & \vdots  & \vdots  & \vdots  & \ddots   \\
\end{matrix} \right).$$
Since
$$\left[ {{x}^{{{q}^{k}}n+i}} \right]\left( \sum\limits_{n=0}^{{{q}^{k}}-1}{{{b}_{n}}{{x}^{n}}} \right)a\left( {{x}^{{{q}^{k}}}} \right)={{a}_{n}}{{b}_{i}}, \qquad 0\le i<{{q}^{k}},$$
then $\left( {{q}^{k}}n+i,{{q}^{k}}m+j \right)$th element of the matrix ${{\left( a\left( x \right)|b\left( x \right) \right)}_{q,k}}$ is equal to
$${{a}_{n-m}}{{b}_{i-j}}{{\left( \begin{matrix}
   {{q}^{k}}n+i  \\
   {{q}^{k}}m+j  \\
\end{matrix} \right)}_{0,q}}={{a}_{n-m}}{{\left( \begin{matrix}
   n  \\
   m  \\
\end{matrix} \right)}_{0,q}}{{b}_{i-j}}{{\left( \begin{matrix}
   i  \\
   j  \\
\end{matrix} \right)}_{0,q}},\qquad0\le i,j<{{q}^{k}}.$$
Thus, ${{\left( a\left( x \right)|b\left( x \right) \right)}_{q,k}}$ is the block matrix, whose $\left( n,m \right)$th block  is the matrix consisting of ${{q}^{k}}$ first rows of the matrix $\left( b\left( x \right),x|{}_{\left[ 0,q \right]}P \right)$ and multiplied by ${{a}_{n-m}}{{n\choose m}_{0,q}}$. Hence
$${{\left( a\left( x \right)|b\left( x \right) \right)}_{q,k}}{{\left( g\left( x \right)|c\left( x \right) \right)}_{q,k}}={{\left( a\left( x \right)\circ g\left( x \right)|b\left( x \right)\circ c\left( x \right) \right)}_{q,k}},$$
$$\left[ {{x}^{n}} \right]a\left( x \right)\circ g\left( x \right)=\sum\limits_{m=0}^{n}{{{\left( \begin{matrix}
   n  \\
   m  \\
\end{matrix} \right)}_{0,q}}{{a}_{m}}{{g}_{n-m}}}.$$
If
$$\left( a\left( x \right),x|{}_{\left[ 0,q \right]}P \right)={{\left( a\left( x \right)|a\left( x \right) \right)}_{q,1}},$$
as in the case of the matrix $_{\left[ 0,q \right]}P$, then $\left( a\left( x \right),x|{}_{\left[ 0,q \right]}P \right)={{\left( a\left( x \right)|a\left( x \right) \right)}_{q,k}}$ for all $k$:
$$a\left( x \right)=\left( \sum\limits_{n=0}^{q-1}{{{a}_{n}}{{x}^{n}}} \right)a\left( {{x}^{q}} \right)=\left( \sum\limits_{n=0}^{{{q}^{k}}-1}{{{a}_{n}}{{x}^{n}}} \right)a\left( {{x}^{{{q}^{k}}}} \right)=\prod\limits_{m=0}^{\infty }{\left( \sum\limits_{n=0}^{{{q}^{k}}-1}{{{a}_{n}}{{x}^{n{{q}^{mk}}}}} \right)},$$
$${{a}_{0}}=1,  \qquad{{a}_{{{q}^{k}}n+i}}={{a}_{n}}{{a}_{i}},  \qquad0\le i<{{q}^{k}},$$ 
$${{a}_{n}}=\prod\limits_{i=0}^{\infty }{{{a}_{{{n}_{i}}}}}, \qquad n=\sum\limits_{i=0}^{\infty }{{{n}_{i}}}{{q}^{i}}={{n}_{0}}+q\left( {{n}_{1}}+q\left( {{n}_{2}}+... \right) \right),  \qquad0\le {{n}_{i}}<q.$$
For example,
$$_{\left[ 0,2 \right]}{{P}^{2}}=\left( \begin{matrix}
   1 & 0 & 0 & 0 & 0 & 0 & 0 & 0 & 0 & 0 & 0 & 0 & 0 & 0 & 0 & 0 & \ldots   \\
   2 & 1 & 0 & 0 & 0 & 0 & 0 & 0 & 0 & 0 & 0 & 0 & 0 & 0 & 0 & 0 & \ldots   \\
   2 & 0 & 1 & 0 & 0 & 0 & 0 & 0 & 0 & 0 & 0 & 0 & 0 & 0 & 0 & 0 & \ldots   \\
   4 & 2 & 2 & 1 & 0 & 0 & 0 & 0 & 0 & 0 & 0 & 0 & 0 & 0 & 0 & 0 & \ldots   \\
   2 & 0 & 0 & 0 & 1 & 0 & 0 & 0 & 0 & 0 & 0 & 0 & 0 & 0 & 0 & 0 & \ldots   \\
   4 & 2 & 0 & 0 & 2 & 1 & 0 & 0 & 0 & 0 & 0 & 0 & 0 & 0 & 0 & 0 & \ldots   \\
   4 & 0 & 2 & 0 & 2 & 0 & 1 & 0 & 0 & 0 & 0 & 0 & 0 & 0 & 0 & 0 & \ldots   \\
   8 & 4 & 4 & 2 & 4 & 2 & 2 & 1 & 0 & 0 & 0 & 0 & 0 & 0 & 0 & 0 & \ldots   \\
   2 & 0 & 0 & 0 & 0 & 0 & 0 & 0 & 1 & 0 & 0 & 0 & 0 & 0 & 0 & 0 & \ldots   \\
   4 & 2 & 0 & 0 & 0 & 0 & 0 & 0 & 2 & 1 & 0 & 0 & 0 & 0 & 0 & 0 & \ldots   \\
   4 & 0 & 2 & 0 & 0 & 0 & 0 & 0 & 2 & 0 & 1 & 0 & 0 & 0 & 0 & 0 & \ldots   \\
   8 & 4 & 4 & 2 & 0 & 0 & 0 & 0 & 4 & 2 & 2 & 1 & 0 & 0 & 0 & 0 & \ldots   \\
   4 & 0 & 0 & 0 & 2 & 0 & 0 & 0 & 2 & 0 & 0 & 0 & 1 & 0 & 0 & 0 & \ldots   \\
   8 & 4 & 0 & 0 & 4 & 2 & 0 & 0 & 4 & 2 & 0 & 0 & 2 & 1 & 0 & 0 & \ldots   \\
   8 & 0 & 4 & 0 & 4 & 0 & 2 & 0 & 4 & 0 & 2 & 0 & 2 & 0 & 1 & 0 & \ldots   \\
   16 & 8 & 8 & 4 & 8 & 4 & 4 & 2 & 8 & 4 & 4 & 2 & 4 & 2 & 2 & 1 & \ldots   \\
   \vdots  & \vdots  & \vdots  & \vdots  & \vdots  & \vdots  & \vdots  & \vdots  & \vdots  & \vdots  & \vdots  & \vdots  & \vdots  & \vdots  & \vdots  & \vdots  & \ddots   \\
\end{matrix} \right).$$
Thus,
$$_{\left[ 0,q \right]}{{P}^{\varphi }}={{\left( {{\left( \frac{1}{1-x} \right)}^{\left( \varphi  \right)}}|{{\left( \frac{1}{1-x} \right)}^{\left( \varphi  \right)}} \right)}_{q,k}},$$
$$\left[ {{x}^{qn+i}} \right]{{\left( \frac{1}{1-x} \right)}^{\left( \varphi  \right)}}=\left( \begin{matrix}
   \varphi +i-1  \\
   i  \\
\end{matrix} \right)\left[ {{x}^{n}} \right]{{\left( \frac{1}{1-x} \right)}^{\left( \varphi  \right)}},  \qquad0\le i<q,$$
$${{\left( _{\left[ 0,q \right]}{{P}^{\varphi }} \right)}_{qn+i,qm+j}}=\left( \begin{matrix}
   \varphi +i-j-1  \\
   i-j  \\
\end{matrix} \right){{\left( _{\left[ 0,q \right]}{{P}^{\varphi }} \right)}_{n,m}}, \qquad0\le i,j<q,$$
$$\left[ {{x}^{n}} \right]{{\left( \frac{1}{1-x} \right)}^{\left( \varphi  \right)}}=\prod\limits_{i=0}^{\infty }{\left( \begin{matrix}
   \varphi +{{n}_{i}}-1  \\
   {{n}_{i}}  \\
\end{matrix} \right)},  \qquad{{\left( _{\left[ 0,q \right]}{{P}^{\varphi }} \right)}_{n,m}}=\prod\limits_{i=0}^{\infty }{\left( \begin{matrix}
   \varphi +{{n}_{i}}-{{m}_{i}}-1  \\
   {{n}_{i}}-{{m}_{i}}  \\
\end{matrix} \right)},$$
$$n=\sum\limits_{i=0}^{\infty }{{{n}_{i}}}{{q}^{i}},  \qquad m=\sum\limits_{i=0}^{\infty }{{{m}_{i}}}{{q}^{i}}  \qquad0\le {{n}_{i}},{{m}_{i}}<q.$$
Denote
$${{\left[ n,\to  \right]}_{\left[ 0,q \right]}}{{P}^{\varphi }}={{u}_{n}}\left( x \right).$$
Then
$${{u}_{0}}\left( x \right)=1,  \qquad{{u}_{qn+m}}\left( x \right)={{w}_{m}}\left( \varphi ,x \right){{u}_{n}}\left( {{x}^{q}} \right), \qquad0\le m<q,$$
$${{u}_{n}}\left( x \right)=\prod\limits_{i=0}^{\infty }{{{w}_{{{n}_{i}}}}}\left( \varphi ,{{x}^{{{q}^{i}}}} \right),  \qquad n=\sum\limits_{i=0}^{\infty }{{{n}_{i}}}{{q}^{i}},  \qquad0\le {{n}_{i}}<q.$$
In particular
$$\left[ n,\to  \right]{}_{\left[ 0,q \right]}{{P}^{-1}}=\prod\limits_{i=0}^{\infty }{{{w}_{{{n}_{i}}}}}\left( -1,{{x}^{{{q}^{i}}}} \right),$$
 $${{w}_{0}}\left( -1,x \right)=1,   \qquad{{w}_{n}}\left( -1,x \right)={{x}^{n-1}}\left( x-1 \right).$$
For exampl,
$${}_{\left[ 0,2 \right]}{{P}^{-1}}=\left( \begin{matrix}
   1 & 0 & 0 & 0 & 0 & 0 & 0 & 0 & 0 & \text{ }0 & \text{   }0 & \ldots   \\
   -1 & 1 & 0 & 0 & 0 & 0 & 0 & 0 & 0 & \text{ }0 & \text{   }0 & \ldots   \\
   -1 & 0 & 1 & 0 & 0 & 0 & 0 & 0 & 0 & \text{ }0 & \text{   }0 & \ldots   \\
   1 & -1 & -1 & 1 & 0 & 0 & 0 & 0 & 0 & \text{ }0 & \text{   }0 & \ldots   \\
   -1 & 0 & 0 & 0 & 1 & 0 & 0 & 0 & 0 & \text{ }0 & \text{   }0 & \ldots   \\
   1 & -1 & 0 & 0 & -1 & 1 & 0 & 0 & 0 & \text{ }0 & \text{   }0 & \ldots   \\
   1 & 0 & -1 & 0 & -1 & 0 & 1 & 0 & 0 & \text{ }0 & \text{   }0 & \ldots   \\
   -1 & 1 & 1 & -1 & 1 & -1\text{ } & -1 & \text{ }1\text{ } & 0 & \text{ }0 & \text{   }0 & \ldots   \\
   -1 & 0 & 0 & 0 & 0 & 0 & 0 & \text{ }0\text{ } & 1 & \text{ }0 & \text{   }0 & \ldots   \\
   1 & -1 & 0 & 0 & 0 & 0 & 0 & \text{ }0\text{ } & -1 & \text{ }1 & \text{   }0 & \ldots   \\
   1 & 0 & -1 & 0 & 0 & 0 & 0 & \text{ }0\text{ } & -1 & \text{ }0 & \text{   }1 & \ldots   \\
   \vdots  & \vdots  & \vdots  & \vdots  & \vdots  & \vdots  & \vdots  & \vdots  & \vdots  & \vdots  & \vdots  & \ddots   \\
\end{matrix} \right).$$
Since
$$\left[ n,\to  \right]{}_{\left[ 0,q \right]}{{P}^{-1}}=\left[ n,\to  \right]{}_{0,q}{{P}^{-1}}, \qquad n<{{q}^{2}},$$
$${{u}_{{{q}^{k}}}}\left( x \right)={{x}^{{{q}^{k}}}}-1,$$
and in other cases polynomials ${{u}_{n}}\left( x \right)$ comprise more than one factor of the form ${{x}^{{{m}_{i}}}}-1$, then
$$_{\left[ 0,q \right]}{{P}^{-1}}x{{\left( 1-x \right)}^{-2}}=\sum\limits_{m=1}^{q-1}{{{x}^{m}}}+\sum\limits_{m=1}^{q-1}{q{{x}^{qm}}}+\sum\limits_{m=2}^{\infty }{{{q}^{m}}{{x}^{{{q}^{m}}}}},$$
$$\log \circ {{\left( 1-x \right)}^{-1}}=\sum\limits_{m=1}^{q-1}{\frac{{{x}^{m}}}{m}}+\sum\limits_{m=1}^{q-1}{\frac{{{x}^{qm}}}{m}}+\sum\limits_{m=2}^{\infty }{{{x}^{{{q}^{m}}}}}.$$
{\bfseries Remark.} Matrices $_{\left[ 0,q \right]}P$, denoted by ${{S}_{q}}$ and called generalized Sierpinski matrices, were introduced in [8]. In the works [8] – [12] properties of these matrices and associated algebras are studied from point of view of the combinatorics and the number theory. In [11] are introduced zero generalized Pascal matrices of the form 
$${{T}^{\left( q \right)}}={{P}_{c\left( x \right)}}\times {}_{\left[ 0,q \right]}P,  \qquad c\left( x \right)=\left( \sum\limits_{n=0}^{q-1}{{{c}_{n}}{{x}^{n}}} \right)c\left( {{x}^{q}} \right), \qquad{{c}_{n}}={{\left( n! \right)}^{-1}}, \qquad0\le n<q,$$
$${{\left( {{T}^{\left( q \right)}} \right)}_{n,m}}={{\left( \begin{matrix}
   n  \\
   m  \\
\end{matrix} \right)}_{q}}=\prod\limits_{i=0}^{\infty }{\left( \begin{matrix}
   {{n}_{i}}  \\
   {{m}_{i}}  \\
\end{matrix} \right)},  \qquad n=\sum\limits_{i=0}^{\infty }{{{n}_{i}}}{{q}^{i}},  \qquad m=\sum\limits_{i=0}^{\infty }{{{m}_{i}}}{{q}^{i}},  \qquad0\le {{n}_{i}},{{m}_{i}}<q,$$
$$\sum\limits_{m=0}^{n}{{{\left( \begin{matrix}
   n  \\
   m  \\
\end{matrix} \right)}_{q}}}{{x}^{m}}=\prod\limits_{i=0}^{\infty }{{{\left( 1+{{x}^{{{q}^{i}}}} \right)}^{{{n}_{i}}}}}.$$

\section{Zero generalized  Riordan group}

Denote
$$\left( a\left( x \right),x|{}_{0}P \right)={{\left( a\left( x \right),1 \right)}_{0}},$$
where $_{0}P$ is a zero generalized Pascal matrix, a particular form of which is specified separately . We construct the matrix  ${{\left( 1,a\left( x \right) \right)}_{0}}$ by the rule
$$[\uparrow ,n]{{\left( 1,a\left( x \right) \right)}_{0}}=[\uparrow ,n]{{\left( {{a}^{\left( n \right)}}\left( x \right),1 \right)}_{0}}={{x}^{n}}\circ {{a}^{\left( n \right)}}\left( x \right).$$
Denote
$${{\left( 1,a\left( x \right) \right)}_{0}}b\left( x \right)=b\circ \left( a\left( x \right) \right),$$
$${{\left( b\left( x \right),1 \right)}_{0}}{{\left( 1,a\left( x \right) \right)}_{0}}={{\left( b\left( x \right),a\left( x \right) \right)}_{0}}.$$
{\bfseries Theorem.}\emph{ Matrices ${{\left( b\left( x \right),a\left( x \right) \right)}_{0}}$, ${{b}_{0}}\ne 0$, ${{a}_{0}}\ne 0$, form a group whose elements are multiplied by the rule}
$${{\left( b\left( x \right),a\left( x \right) \right)}_{0}}{{\left( f\left( x \right),g\left( x \right) \right)}_{0}}={{\left( b\left( x \right)\circ f\circ \left( a\left( x \right) \right),a\left( x \right)\circ g\circ \left( a\left( x \right) \right) \right)}_{0}}.$$
{\bfseries Proof.} Let ${{\left( _{0}P \right)}_{n,m}}=\left( n,m \right)$. Then
$${{x}^{m}}\circ {{x}^{n}}=\left( m+n,n \right){{x}^{m+n}} ,  \qquad{{x}^{m}}\circ b\left( x \right)=\sum\limits_{n=0}^{\infty }{{{b}_{n}}}\left( m+n,n \right){{x}^{m+n}},$$
$${{\left( 1,a\left( x \right) \right)}_{0}}{{x}^{m}}\circ b\left( x \right)=\sum\limits_{n=0}^{\infty }{{{b}_{n}}}\left( m+n,n \right){{x}^{m+n}}\circ {{a}^{\left( m+n \right)}}\left( x \right)=$$
$$={{x}^{m}}\circ {{a}^{\left( m \right)}}\left( x \right)\circ \sum\limits_{n=0}^{\infty }{{{b}_{n}}}{{x}^{n}}\circ {{a}^{\left( n \right)}}\left( x \right)={{x}^{m}}\circ {{a}^{\left( m \right)}}\left( x \right)\circ b\circ \left( a\left( x \right) \right),$$
or
$${{\left( 1,a\left( x \right) \right)}_{0}}{{\left( b\left( x \right),1 \right)}_{0}}={{\left( b\circ \left( a\left( x \right) \right),a\left( x \right) \right)}_{0}}.$$
Then
$${{\left( 1,a\left( x \right) \right)}_{0}}b\left( x \right)\circ c\left( x \right)=b\circ \left( a\left( x \right) \right)\circ c\circ \left( a\left( x \right) \right),$$
$${{\left( 1,a\left( x \right) \right)}_{0}}{{x}^{m}}\circ {{b}^{\left( m \right)}}\left( x \right)={{x}^{m}}\circ {{a}^{\left( m \right)}}\left( x \right)\circ {{\left( b\circ\left( a\left( x \right) \right) \right)}^{\left( m \right)}},$$
or
$${{\left( 1,a\left( x \right) \right)}_{0}}{{\left( 1,b\left( x \right) \right)}_{0}}={{\left( 1,a\left( x \right)\circ b\circ \left( a\left( x \right) \right) \right)}_{0}}.$$ 
{\bfseries Remark 1.} If $_{0}P={}_{0,q}P\times {{P}_{c\left( x \right)}}$, matrices
$${{\left( b\left( x \right),a\left( x \right) \right)}_{0}}, \qquad{{b}_{qn+m}}=0, \qquad{{a}_{qn+m}}=0, \qquad0<m<q,$$
form a subgroup common to all groups associated with the set of generalized Pascal matrices $_{\varphi ,q}P\times {{P}_{c\left( x \right)}}$. For example, if $_{0}P={}_{0,2}P$,
$$\left( \begin{matrix}
   1 & 0 & 0 & 0 & 0 & 0 & \ldots   \\
   0 & 1 & 0 & 0 & 0 & 0 &  \ldots   \\
   0 & 0 & 1 & 0 & 0 & 0 & \ldots   \\
   0 & 1 & 0 & 1 & 0 & 0 &  \ldots   \\
   0 & 0 & 2 & 0 & 1 & 0 &  \ldots   \\
   0 & 1 & 0 & 3 & 0 & 1 & \ldots   \\
   \vdots  & \vdots  & \vdots  & \vdots  & \vdots  & \vdots  & \ddots   \\
\end{matrix} \right)\left( \begin{matrix}
   1 & 0 & 0 & 0 & 0 & 0 &  \ldots   \\
   0 & 1 & 0 & 0 & 0 & 0 &  \ldots   \\
   0 & 0 & 1 & 0 & 0 & 0 & \ldots   \\
   0 & 1 & 0 & 1 & 0 & 0 &  \ldots   \\
   0 & 0 & 2 & 0 & 1 & 0 & \ldots   \\
   0 & 1 & 0 & 3 & 0 & 1 & \ldots   \\
   \vdots  & \vdots  & \vdots  & \vdots  & \vdots  & \vdots  & \ddots   \\
\end{matrix} \right)=\left( \begin{matrix}
   1 & 0 & 0 & 0 & 0 & 0 &  \ldots   \\
   0 & 1 & 0 & 0 & 0 & 0 &  \ldots   \\
   0 & 0 & 1 & 0 & 0 & 0 & \ldots   \\
   0 & 2 & 0 & 1 & 0 & 0 &  \ldots   \\
   0 & 0 & 4 & 0 & 1 & 0 &  \ldots   \\
   0 & 5 & 0 & 6 & 0 & 1 &  \ldots   \\
   \vdots  & \vdots  & \vdots  & \vdots  & \vdots  & \vdots  &  \ddots   \\
\end{matrix} \right),$$
or
$${{\left( 1,\frac{1}{1-{{x}^{2}}} \right)}_{0}}{{\left( 1,\frac{1}{1-{{x}^{2}}} \right)}_{0}}={{\left( 1,\frac{1-{{x}^{2}}}{1-3{{x}^{2}}+{{x}^{4}}} \right)}_{0}}=\left( 1,\frac{x\left( 1-{{x}^{2}} \right)}{1-3{{x}^{2}}+{{x}^{4}}} \right),$$
where the all  matrices are ordinary Riordan arrays. For comparison,
$$\left( \begin{matrix}
   1 & 0 & 0 & 0 & 0 & 0 &  \ldots   \\
   0 & 1 & 0 & 0 & 0 & 0 &  \ldots   \\
   0 & 0 & 1 & 0 & 0 & 0 &  \ldots   \\
   0 & 1 & 0 & 1 & 0 & 0 & \ldots   \\
   0 & 0 & 2 & 0 & 1 & 0 &  \ldots   \\
   0 & 1 & 0 & 3 & 0 & 1 & \ldots   \\
   \vdots  & \vdots  & \vdots  & \vdots  & \vdots  & \vdots  & \ddots   \\
\end{matrix} \right)\left( \begin{matrix}
   1 & 0 & 0 & 0 & 0 & 0 &  \ldots   \\
   0 & 1 & 0 & 0 & 0 & 0 &  \ldots   \\
   0 & 0 & 1 & 0 & 0 & 0 &  \ldots   \\
   0 & 1 & 2 & 1 & 0 & 0 & \ldots   \\
   0 & 0 & 2 & 0 & 1 & 0 & \ldots   \\
   0 & 1 & 4 & 3 & 4 & 1 & \ldots   \\
   \vdots  & \vdots  & \vdots  & \vdots  & \vdots  & \vdots  &  \ddots   \\
\end{matrix} \right)=\left( \begin{matrix}
   1 & 0 & 0 & 0 & 0 & 0 &  \ldots   \\
   0 & 1 & 0 & 0 & 0 & 0 &  \ldots   \\
   0 & 0 & 1 & 0 & 0 & 0 &  \ldots   \\
   0 & 2 & 2 & 1 & 0 & 0 &  \ldots   \\
   0 & 0 & 4 & 0 & 1 & 0 & \ldots   \\
   0 & 5 & 10 & 6 & 4 & 1 &  \ldots   \\
   \vdots  & \vdots  & \vdots  & \vdots  & \vdots  & \vdots  &  \ddots   \\
\end{matrix} \right),$$
or
$${{\left( 1,\frac{1}{1-{{x}^{2}}} \right)}_{0}}{{\left( 1,\frac{1}{1-x} \right)}_{0}}={{\left( 1,\frac{1}{1-x-{{x}^{2}}} \right)}_{0}},$$
where
$$a\left( x \right)=\frac{1}{1-{{x}^{2}}},   \qquad b\left( x \right)=\frac{1}{1-x},$$
$$b\circ \left( a\left( x \right) \right)=\frac{1-{{x}^{2}}}{1-x-{{x}^{2}}} ,   \qquad a\left( x \right)\circ b\circ \left( a\left( x \right) \right)=\frac{1}{1-x-{{x}^{2}}}.$$
{\bfseries Remark 2.} If in the algebra associated with the matrix $_{0}P$, the series ${{a}_{1}}\left( x \right)$, ${{a}_{2}}\left( x \right)$ belong to the same group of $l$-series, i.e. have the form (5), then
$$[\uparrow ,n]{{\left( 1,{{a}_{1}}\left( x \right) \right)}_{0}}={{x}^{n}},\qquad\mbox {if} \qquad  \left[ {{x}^{n}} \right]{{a}_{2}}\left( x \right)\ne 0,$$
$${{a}_{2}}\circ \left( {{a}_{1}}\left( x \right) \right)={{a}_{2}}\left( x \right),  \qquad{{\left( 1,{{a}_{1}}\left( x \right) \right)}_{0}}{{\left( 1,{{a}_{2}}\left( x \right) \right)}_{0}}={{\left( 1,{{a}_{1}}\left( x \right)\circ {{a}_{2}}\left( x \right) \right)}_{0}}.$$
{\bfseries Remark 3.} Identities (3), (4) provide following analogue of  the Lagrange inversion theorem:\\
{\bfseries Theorem.}\emph{ If the matrices ${{\left( 1,{{a}^{\left( -1 \right)}}\left( x \right) \right)}_{0}}$, ${{\left( 1,b\left( x \right) \right)}_{0}}$ are mutually inverse, then }
$$\left[ {{x}^{n}} \right]{{b}^{\left( \varphi  \right)}}\left( x \right)=\frac{\varphi }{\varphi +n}\left[ {{x}^{n}} \right]{{a}^{\left( \varphi +n \right)}}\left( x \right).$$
{\bfseries Proof.} Denote
$$\left[ {{x}^{n}} \right]\left( 1-x{{\left( \log \circ a\left( x \right) \right)}^{\prime }} \right)\circ {{a}^{\left( m \right)}}\left( x \right)=c_{n}^{m}, \qquad{{c}_{m}}\left( x \right)=\sum\limits_{n=0}^{\infty }{c_{n}^{m+n}{{x}^{n}}}.$$
Construct the matrix $C$,  $\left[ \uparrow ,n \right]C={{x}^{n}}\circ {{c}_{n}}\left( x \right)$:
$$C=\left( \begin{matrix}
   c_{0}^{0} & 0 & 0 & 0 & \ldots   \\
   c_{1}^{1} & c_{0}^{1} & 0 & 0 & \ldots   \\
   c_{2}^{2} & c_{1}^{2} & c_{0}^{2} & 0 & \ldots   \\
   c_{3}^{3} & c_{2}^{3} & c_{1}^{3} & c_{0}^{3} & \ldots   \\
   \vdots  & \vdots  & \vdots  & \vdots  & \ddots   \\
\end{matrix} \right)\times {}_{0}P.$$
It's obvious that
$$\left[ n,\to  \right]C=\left[ n,\to  \right]{{\left( \left( 1-x{{\left( \log \circ a\left( x \right) \right)}^{\prime }} \right)\circ {{a}^{\left( n \right)}}\left( x \right),1 \right)}_{0}}.$$
Since
$$\left( 1-x{a}'\left( x \right)\circ {{a}^{\left( -1 \right)}}\left( x \right) \right)\circ {{a}^{\left( m \right)}}\left( x \right)={{a}^{\left( m \right)}}\left( x \right)-\frac{x}{m}{{\left( {{a}^{\left( m \right)}}\left( x \right) \right)}^{\prime }},$$
or
$$\left[ {{x}^{n}} \right]\left( 1-x{{\left( \log \circ a\left( x \right) \right)}^{\prime }} \right)\circ {{a}^{\left( m \right)}}\left( x \right)=\frac{m-n}{m}\left[ {{x}^{n}} \right]{{a}^{\left( m \right)}}\left( x \right),$$
then
$$\left[ {{x}^{m+n}} \right]C{{x}^{m}}\circ {{a}^{\left( -m \right)}}\left( x \right)=\left[ {{x}^{n}} \right]\left( 1-x{{\left( \log \circ a\left( x \right) \right)}^{\prime }} \right)\circ {{a}^{\left( n \right)}}\left( x \right)=1,\quad n=0;\quad =0,\quad n>0.$$
Thus,
$$C={{\left( 1,b\left( x \right) \right)}_{0}},  \qquad\left[ {{x}^{n}} \right]{{b}^{\left( m \right)}}\left( x \right)=c_{n}^{m+n}=\frac{m}{m+n}\left[ {{x}^{n}} \right]{{a}^{\left( m+n \right)}}\left( x \right).$$
Since
$${{b}^{\left( \varphi  \right)}}\left( x \right)=\sum\limits_{n=0}^{\infty }{\frac{{{\varphi }^{n}}}{n!}}{{\left( \log \circ b\left( x \right) \right)}^{\left( n \right)}}=\sum\limits_{n=0}^{\infty }{{{b}_{n}}}\left( \varphi  \right){{x}^{n}},$$

$${{a}^{\left( \varphi  \right)}}\left( x \right)=\sum\limits_{n=0}^{\infty }{\frac{{{\varphi }^{n}}}{n!}}{{\left( \log \circ a\left( x \right) \right)}^{\left( n \right)}}=\sum\limits_{n=0}^{\infty }{{{a}_{n}}}\left( \varphi  \right){{x}^{n}},$$
where ${{b}_{n}}\left( x \right)$, ${{a}_{n}}\left( x \right)$ are polynomials, then

$${{b}_{n}}\left( x \right)=x{{\left( x+n \right)}^{-1}}{{a}_{n}}\left( x+n \right).$$

E-mail: {evgeniy\symbol{"5F}burlachenko@list.ru}
\end{document}